\documentclass[nomath,noamsfonts]{amsart}
\input{psfig.tex}
\input{psfig.sty}
\begin{document}
\title{Cubulations, immersions, mappability and a problem of Habegger}
\author{Louis Funar}
\address{Institut Fourier, BP 74, Univ.Grenoble I, Math\'ematiques,  
38402 Saint-Martin-d'H\`eres cedex, France}  
\email{funar@fourier.ujf-grenoble.fr}
\date{\today}
\subjclass{Primary  57 Q 15, 57 R 42, $\;$  Secondary 51 M 20, 57 Q 35}. 
\begin{abstract}
The aim of this paper (inspired from a problem of Habegger) 
is to describe the set of cubical decompositions
of  compact manifolds mod out by  a set of combinatorial moves analogous to 
the bistellar moves considered by Pachner,  which we call bubble
moves. One constructs a surjection from this set onto the 
the bordism group of codimension one immersions in the manifold. 
The connected sums of manifolds and immersions induce multiplicative 
structures which are respected by this surjection. 
We prove that  those cubulations  which map  combinatorially into the 
standard decomposition of ${\bf R}^n$ for large enough $n$ 
(called mappable), are equivalent.  
Finally we classify the cubulations of the 2-sphere. 
\end{abstract}
\maketitle

\section{Introduction and statement of results}
\subsection{Outline}
Stellar moves were first considered by Alexander  (\cite{Al}) who 
proved that they can relate any two  triangulations of a polyhedron. 
Alexander's moves were refined to a finite set of
local (bistellar) moves  which still 
act transitively  on the triangulations of manifolds,
according to  Pachner (\cite{pachner}). 
Using Pachner's result Turaev and Viro 
proved that certain state-sums associated to a triangulation yield
topological invariants of 3-manifolds (see \cite{Tu}).  
Recall that a bistellar move 
(in dimension $n$) excises $B$ and replaces it with 
$B'$, where $B$ and $B'$ are complementary balls, subcomplexes 
in the boundary of the standard $(n+1)$-simplex.
For a nice exposition of Pachner's result and various extensions, see 
\cite{Li}.

The Turaev-Viro invariants
carry less information than the Reshetikhin-Turaev 
invariants, which are defined using  Dehn surgery presentations 
instead of triangulations. 
In fact the latter have a strong 4-dimensional flavor,
as explained by the theory of shadows developed by Turaev (see
\cite{Tu}). This motivates the study of 
state-sums based on cubulations, as an alternative  way to get
intrinsic invariants possibly containing more 
information (e.g. the phase factor).  
A {\it cubical complex} is a complex  $K$ consisting of Euclidean cubes, 
such that the intersection of two cubes  is a finite union of cubes
from  $K$, once a cube is in $K$ all its faces are still in $K$, and no
identifications of faces of the same cube are allowed. 
A {\it cubulation} of a manifold is specified by a cubical
complex PL homeomorphic to the manifold.   
In order to apply the state-sum machinery to these 
decompositions we need  an analogue of 
Pachner's theorem. Specifically, N.Habegger asked (see    
problem 5.13 from R.Kirby's list (\cite{kirby})) the following:
\newtheorem{alba}{Problem}
\begin{alba}
Suppose $M$ and $N$ are  PL-homeomorphic cubulated $n$-manifolds. Are
they related by the following set of moves: excise $B$ and replace it
by
 $B'$, 
where $B$ and $B'$ are complementary balls (union of $n$-cubes) 
in the boundary of the standard $(n+1)$-cube?
\end{alba}
These moves will be called {\em bubble} moves in the sequel.
Among them, those for which  $B$ or $B'$ does not contain parallel 
(when viewed in the $n+1$-cube) faces are called {\em np-bubble} moves. 
There are $n+1$ distinct np-bubble moves $b_k$, $k=1,2,...,n+1$ and their 
inverses, where the support $B$ of $b_k$ is the union of exactly $k$ cubes.
For $n=2$ there is one bubble move which is not a np-bubble (see picture 1). 
Set $C(M)$ for  the  set of  cubulations of a closed manifold $M$,  
$CBB(M)$ for the 
equivalence classes of cubulations mod np-bubble moves and  $CB(M)$ for the 
equivalence classes of cubulations mod bubble moves. 
The answer to Habegger's question, as it states, 
is negative because the triangle and the square are not bubble equivalent. 
In fact, for $n=1$ the move $b_1$ divides an edge into three edges 
and so $CB(S^1) =  CBB(S^1)={\bf Z}/2{\bf Z}$.  
Therefore a complete answer would rather consist of a description of $CB(M)$. 
Another way is to avoid the difficulties of a direct approach 
by looking for a sufficiently large
class of cubulations having an intrinsic characterization and within
which the cubulations are  equivalent. 
The aim of this paper is to formulate some partial solutions 
along these lines.  

For instance, one associates  to each cubulation $C$ of $M$, a codimension 
one normal crossings  immersion $\varphi_C$ in $M$. 
In this way one obtains a surjective map from the set of  
(marked) cubulations mod bubble moves to the bordism set of  
immersions. The latter has a homotopical description via the 
Pontryagin-Thom construction. We conjecture that this surjection 
is a bijection.
On the other hand let us restrict to cubulations which can be combinatorially
mapped into the standard cubulation of some Euclidean space (called 
mappable). One can approximate an
ambient isotopy between two cubulations by some cubical 
sub-complexes of the standard cubulation.  
The path of cubical approximations  is 
locally constant except for a finite number of critical values of the 
parameter, when a jump described by a bubble move occurs. 
As a consequence two  mappable  cubulations are 
bubble equivalent.   
We prove that the connected sum of cubulations mod bubble moves is 
well-defined, and this is compatible with the composition map for 
immersions. Finally we consider the case of $CB(S^2)$ and show by a
direct combinatorial proof that $CB(S^2)={\bf Z}/2{\bf Z}$. 
\begin{figure}
\hspace{0.5cm}\psfig{figure=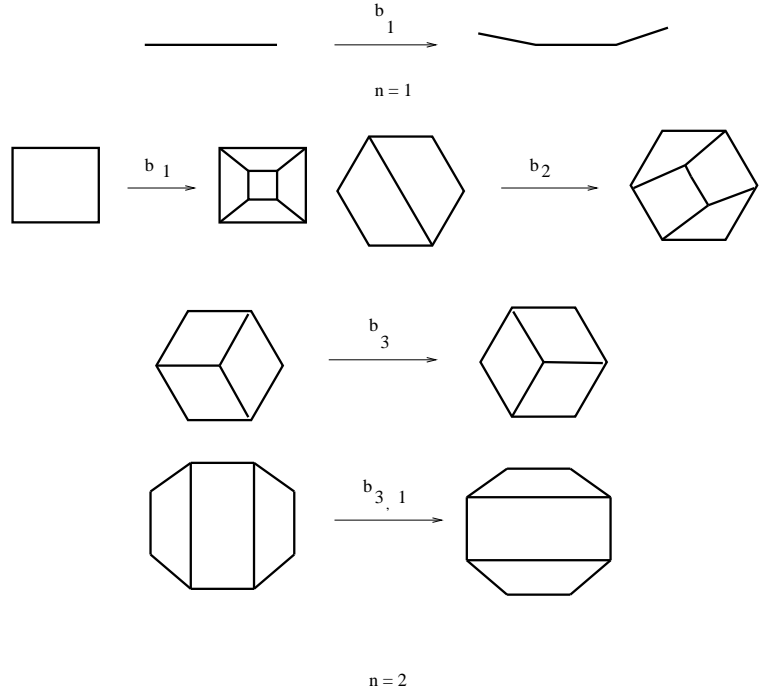,width=10cm}
\caption{Bubble moves for $n=1$ and $n=2$}
\end{figure}

{\bf Acknowledgements}: Part of this work was done during the author's
visit at University of Palermo and University of Columbia, whose
support and hospitality are gratefully acknowledged. I'm thankful to 
E.Babson, J.Birman, C.Blanchet, R.Casali, C.Chan,
L.Guillou, N.Habegger, T.Kashiwabara, A.Marin, D.Matei, S.Matveev, 
V.Po\'enaru, R.Popescu, V.Sergiescu and the referee for  helpful
discussions, suggestions and improvements. 

\subsection{Elementary obstructions}
We outline the combinatorial approach in higher dimensions from 
\cite{Fun}. For a cubulation $x\in C(M)$ of the 
$n$-manifold $M$  the component $f_i(x)$ of the $f$-vector ${\bf f}(x)$ 
counts the number of $i$-dimensional cubes in $x$. 
The orbit of the $f$-vector ${\bf f}$ under bubble moves has the form 
${\bf f}+ \Lambda(n) \subset {\bf  Z}^{n+1}$, 
where $\Lambda(n)$ is a lattice. Therefore we have 
an induced map $CB(M)\to {\bf  Z}^{n+1}/\Lambda(n)$ taking values in 
a finite Abelian group. 

\newtheorem{alba1}{Proposition}[section]
\begin{alba1}
There exist nonzero even numbers $a_i(n)\in {\bf Z}_+$ such that
the projection ${\bf
  Z}^{n+1}/\Lambda(n)\longrightarrow \prod_{i=0}^n {\bf Z}/a_i(n){\bf
  Z}$
is surjective.  The greatest such numbers $a_i(n)$ verify 
$a_n(n)=2$, $a_{n-1}(n)=2n$, $a_{n-2}(n)=2$, $a_0(n)=2$, $a_1(n)=
3+(-1)^n, (n>2)$. 
\end{alba1}
See \cite{Fun} for the proof. 
Let $fb$ be the class of ${\bf f}$ in 
$\prod_{i=0}^n {\bf Z}/a_i(n){\bf
  Z}$  and  $fb^{(2)}$ be the  reduced elements modulo
$(2,2,2,...,2,2n,2)$. Notice that $\Lambda(n)$ is not a product
lattice in general.  For instance, when  $n=3$ there is an 
additional invariant $f_0+f_1\in {\bf Z}/4{\bf Z}$. 

A  natural problem is  to compute the image 
$fb(CB(M))$ for given $M$. Some partial
results for the  mod 2 reductions 
$fb^{(2)}(CB(M))$ are known. 
This is equivalent to characterize those $f$-vectors mod 2
which can be realized by cubulations of the manifold $M$.  
There are constraints for the existence of a simplicial 
polyhedron with a given $f$-vector and fixed topological type. 
For convex simplicial polytopes one has  McMullen's conditions (see 
\cite{mcm,barnette, bill-lee, stanley,barnette1,mcm2}). 
The complete characterization of the $f$-vectors of simplicial
polytopes (and PL-spheres) was obtained in \cite{stanley2}. 
The analogous problem of the realization of $f$-vectors by cubical
polytopes  has also been addressed in some recent papers, for example
\cite{BB,BC,He,J} and references therein.   
The new feature is that, unlike for the simplicial case, there are
parity restrictions on the $f$-vectors (see 
\cite{BB}).  
The relationship between cubical PL $n$-spheres and
the immersions was described in the following result of  
Babson and Chan (see \cite{BC}):
\newtheorem{albat}[alba1]{Proposition}
\begin{albat}  
There exists a  cubical $n$-sphere $K$ with given 
$f_i(mod \; 2)$  if and
only if there exists a codimension 1 normal crossings immersion
$\varphi: M\longrightarrow S^n$ such that 
$f_i(K)= \chi (X_i(M,\varphi)) (mod \;\; 2)$,
where $\chi$ denotes the Euler characteristics, and 
$X_i(M,\varphi)$ is the set of $i$-tuple points.
\end{albat}
There is a wide literature on  immersions, and especially on the
function $\theta_n$ counting the number of multiple $n$-points mod 2, 
which was
considered first  by Freedman (\cite{F}). 
Earlier Banchoff \cite{B} has proved 
that the number of triple points of a closed surface $S$ immersed
in ${\bf R}^3$ is $\chi(S)(mod\; 2)$. There
is an induced  
homomorphism  $\theta_n : B_n \longrightarrow {\bf Z}/2{\bf Z}$ on
the Abelian group $B_n$ of bordism classes of immersions of $(n-1)$-manifolds 
in $S^{n}$. Now  $\theta_n$ is surjective (i.e. nontrivial) if and
only if $fb^{(2)}_{n-1}(S^n)={\bf Z}/2{\bf Z}$. From the 
results concerning the function $\theta_n$ obtained in 
\cite{F,E1,E2,E3,H,Kos,KS,L,carter1,carter2,carter3} we deduce that
the $f$-vectors of a $n$-sphere   
have the following properties (see also \cite{BC}):
\begin{enumerate}
\item For $n=2$ we have $f_0=f_2(mod \; 2)$ and $f_1=0(mod \; 2)$ and 
thus 
$fb^{(2)}(CB(S^2))=fb^{(2)}(CBB(S^2))={\bf Z}/2{\bf Z}$.
\item For $n=3$, $f_0=f_1=0(mod \; 2)$, $f_2=f_3(mod \; 2)$. 
The Boy  immersion $j:{\bf R}P^2\longrightarrow S^3$  
has a single triple point and so there exists
a PL 3-sphere with an odd number of facets. Therefore 
$fb^{(2)}(CB(S^3))=fb^{(2)}(CBB(S^3))={\bf Z}/2{\bf Z}$. 
\item The characterization of  $fb^{(2)}_{n-1}(S^n)$ is
reduced to a homotopy problem: $fb^{(2)}_{n-1}(S^n)={\bf Z}/2{\bf Z}$ 
if and only if 
\begin{enumerate}
\item either $n$ is $1,3,4$ or $7$.
\item or else $n=2^a-2$, with $a\in {\bf Z}_+$, and there exists a  framed 
$n$-manifold with Kervaire invariant 1. The latter is known to be true
for $n=2,6,14,30,62$.
\end{enumerate}
\item If we consider only  edge-orientable cubulations 
(see \cite{He}) then  $fb^{(2)}_{n-1}(S^n)$ is known.
The edge-orientability is equivalent to the orientability of the 
manifold immersed in $S^n$ and the
restriction of the map $\theta_n$ to the subgroup of oriented bordism
classes was computed in \cite{F}. 
In particular $f_{n-1}=0(mod \; 2)$ if $n\neq 1,2,4$. 
\end{enumerate}
\subsection{The 2-dimensional case}
To a surface cubulation we can associate a set of immersed 
circles $K_i$ obtained from the 
union of arcs joining the opposite sides in each square.
The cubulation is {\em simple} if the circles
$K_i$ are individually embedded in the respective surface. 
Simple is equivalent to mappable for the cubulations of  $S^2$ (see
below). A cubulation is called 
{\em semi-simple} if 
each image  circle $\varphi(K_i)$ has an even number of double points,
which form cancelling pairs. Two double points form  a {\em cancelling
  pair} if they are connected by two distinct and disjoint arcs.
\newtheorem{alba11}[alba1]{Theorem}
\begin{alba11}
The np-bubble moves act transitively on the set of simple 
cubulations of  $S^2$. The orbit of the standard  cubulation 
is the set of semi-simple cubulations. 
The map $fb^{(2)}=f_0(mod \; 2)$ is an isomorphism between 
$CB(S^2)$ and ${\bf Z}/2{\bf Z}$. 
\end{alba11}

\subsection{Bordisms of immersions}
Let us consider the set $\mathcal{I} (M)$ (respectively 
$\mathcal{I}^+ (M)$ in the orientable case) of  bordisms of
codimension 1 nc (i.e. normal crossings) immersions in the 
manifold $M$. Two nc-immersions 
$f_i:N_i\longrightarrow M$ of the $(n-1)$-manifolds 
$N_i$ are {\em bordant} if there exists a proper  nc-immersion 
$f:N\longrightarrow M\times [0,1]$ of some 
cobordism $N$ between $N_1$ and $N_2$, such that 
the restriction of $f$ to $N_i$ is isotopic to $f_i$. Using general
position arguments one may get rid of the nc-assumption.

A {\em marked} cubulation is a cubulation $C$ of the 
manifold $M$, endowed with a PL-homeomorphism  $|C|\longrightarrow M$ of
its subjacent space $|C|$, considered up to isotopy.
If a bubble move is performed on $C$, then there is a natural 
marking induced for the bubbled cubulation. Thus it makes sense to 
consider the set  $\widetilde{CB}(M)$ of marked 
cubulations mod bubble moves.  
We  associate to each marked cubulation $C$  a codimension 1 nc-immersion 
$\varphi_C: N_C\longrightarrow M$ (the cubical complex $N_C$ 
was called the derivative complex in \cite{BC}).
Each cube is divided into $2^n$ equal cubes by $n$ hyperplanes 
which we call sections. When gluing together cubes in a cubical
complex the sections are glued accordingly.  The  union of 
the  hyperplane sections  form the image of a 
codimension 1 nc-immersion. In the differentiable case one  uses  
a suitable smoothing when gluing the faces. 
If the cubulation $C$ is edge-orientable (see \cite{He}), and $M$ is oriented then $N_C$
is an oriented manifold.
\newtheorem{alba22}[alba1]{Theorem}
\begin{alba22}
The map $C\to \varphi_C$ induces a surjection
$I:\widetilde{CB}(M)\longrightarrow\mathcal{I}(M)$.
\end{alba22}
The Theorem 1.3 says that the map  $I$ is injective for $M=S^2$. 
We conjecture  that  
$I$  is bijective. In particular 
$\widetilde{CB}(M)$  would depend only on the homotopy type of $M$ and  
the functor $\widetilde{CB}$ which associates to $M$ the set
$\widetilde{CB}(M)$ would be (homotopically) representable.  
Notice that $CB(M)=\widetilde{CB}(M)/\mathcal M (M)$, where 
$\mathcal M (M)$ is the mapping class group of $M$, i.e. the group of 
homeomorphisms of $M$ up to isotopy. Using the classical Pontryagin-Thom 
construction (see e.g. \cite{V}) it follows that 
$\mathcal{I}(M)=[M_c, \Omega^{\infty}S^{\infty}{\bf RP}^{\infty}]$,
and 
$\mathcal{I}^+(M)=[M_c, \Omega^{\infty}S^{\infty}S^1]$, where $M_c$ is
the one point compactification, $\Omega$ denotes the loop space, 
$S$ the reduced suspension and the brackets denote the set of the homotopy
classes of maps. Moreover $\mathcal{I}^+(M)=\pi^{1}(M_c) $ 
is the first cohomotopy group $\pi^{1}(M_c)$. 
The cohomotopy groups of spheres 
can be  computed: 
$\mathcal I(S^n)=\pi_n^s({\bf RP}^{\infty})$, where 
 $\pi_n^s({\bf RP}^{\infty})$ is  the $n$-th stable homotopy group,
 and  $\mathcal I^+(S^n)=\pi_n^s(S^1)=\pi_{n-1}^s$. 
It is known that 
$\pi_1^s({\bf RP}^{\infty})=\pi_2^s({\bf RP}^{\infty})={\bf Z}/2{\bf
  Z}$, 
$\pi_3^s({\bf RP}^{\infty})={\bf Z}/8{\bf Z}$, 
and a few values of the stable stems are tabulated below: 

\vspace{0.5cm}
\begin{center}
\begin{tabular}{|c|c|c|c|c|c|c|c|c|c|}\hline
$n$ & 0 & 1 & 2 & 3 & 4 & 5 & 6 & 7 & 8 \\ \hline
$\pi_{n}^s$ & ${\bf Z}$ &$ {\bf Z}/2$ & ${\bf Z}/2$ & ${\bf Z}/24$
&$ 0$ &$ 0$ & ${\bf Z}/2$ & ${\bf Z}/240$ &$ {\bf Z}/2\oplus {\bf Z}/2 $\\ \hline
\end{tabular}
\end{center}
\vspace{0.5cm}

Let us introduce now the set $\mathcal{C}(M)$ of bordisms of 
cubulations of the manifold
$M$.  The cubulations $C_1$ and $C_2$ are {\em bordant} if there exists 
a cubulation $C$ of $M\times[0,1]$ whose restrictions on the 
boundaries are the $C_i$.  The identity induces a map $CB(M)\longrightarrow
\mathcal{C}(M)$. The question on the existence of  an 
inverse arrow is similar to  Wall's theorem about the
existence of formal deformations between  simple homotopy equivalent 
$n$-complexes through $(n+1)$-complexes (for $n\neq 2$). Remark
that any two cubulations become bordant 
when  suitably subdivided. Consider some  cubulation 
of the sphere $S^n$ which is bubble equivalent to the standard one. 
We can view the bubble moves as the result of gluing and deleting 
$(n+1)$-cubes (after some thickening) to the given cubulation. 
It follows that any such cubulation bounds, i.e. it is the boundary
of a cubulation of the $(n+1)$-ball. For instance  a
polygon bounds iff it has an even number of edges. 
For $n=2, 3$ it might be true  that  the boundary of a ball cubulation is 
bubble equivalent to the standard one, but the result 
cannot be extended to $n\geq 4$. This is analogous  to  the existence of   
non-shellable triangulations of the ball for $n\geq 3$. 
We define a {\em shuffling}  to be a 
sequence of moves where  shellings (i.e. adding iteratively cells,  
each intersecting the  union of the previous ones along a ball) alternate 
with  inverse shellings. An equivalent statement of 
 Pachner's theorem is that all
triangulations can be shuffled.  In the case of  cubulations the 
first  obstruction  for shuffling is that the cubulation bounds. 
However there exist cubulations which bound but cannot be shuffled for 
$n=4$. Consider for example the connected sum  $x\sharp x$, where 
$I(x)$ is the generator of the third stable stem.
We will prove below that the connected sum of cubulations 
makes $CB(S^n)$ a monoid. But $x\sharp x$ bounds and
if it can be shuffled then $x$ should have  order 2. 
This is impossible because $I$ is a monoid homomorphism and  
the bordism group is $\mathcal I^+(S^4)={\bf Z}/24{\bf Z}$. 

\subsection{Embeddable and mappable cubulations}
Cubical complexes, as objects of study from a topological point of
view, were also considered  by 
Novikov (\cite{novikov}, p.42) which asked whether  
a cubical complex of dimension $n$ embeds in (or can be 
mapped to) the $n$-skeleton of the standard cubic lattice of some dimension
$N$. These  are called {\em embeddable} and respectively {\em mappable} 
cubulations. By the  standard cubic lattice 
(or the standard  cubical decomposition) is meant the usual partition
of ${\bf R}^N$ 
into cubes with vertices in  ${\bf Z}^N$.
Several results were obtained in \cite{dss1,dss2,dss3,dss,kara}.
\newtheorem{alba3}[alba1]{Theorem}
\begin{alba3}
The mappable cubulations  of a PL manifold $M$ are
bubble equivalent. 
\end{alba3}
Let us say that a cubulation is {\em simple} if no path in which
consecutive points correspond to edges which are opposite sides of
some square of the cubulation, contains two orthogonal edges 
from the same cube.
The cubulation is {\em standard} if any two of its cubes  are either
disjoint  or have exactly one common face.
An immediate  observation is that embeddable cubulations are standard and
simple and mappable cubulations are simple. 
On the other hand the simplicity is very close to the mappability, at
least for manifolds with small fundamental group. We have for instance
the following  results of Karalashvili (\cite{kara}) and 
Dolbilin, Shtanko and Shtogrin (\cite{dss3}): 
\begin{enumerate}
\item The double (i.e. the result of dividing each $k$-dimensional cube 
in $2^k$ equal cubes)  of a simple cubulation is mappable. 
\item A simple cubulation of a manifold $M$ satisfying 
$H^1(M,{\bf Z}/2{\bf Z})=0$ is mappable. 
\item From a simplicial decomposition $S$ one constructs a cubulation 
$C(S)$ by dividing  each $n$-simplex into $n+1$ cubes. 
Then the cubical decomposition $C(S)$ is embeddable.
\end{enumerate}
In particular  cubulations coming from triangulations  are 
bubble equivalent. Also the simple cubulations of the sphere are
equivalent.  Notice 
that the set of simple (or mappable) cubulations is not closed to
arbitrary  bubble moves. In general the simplicity is not preserved 
by the move $b_2$. 

\subsection{Multiplicative structures}
The connected sum of the manifolds $M$ and $N$ is denoted by $M\sharp
N$. When appropriately extended  to cubulations $\sharp$ depends on  various
choices, but after passing to bubble equivalence classes these ambiguities
disappear. 
\newtheorem{alba5}[alba1]{Theorem}
\begin{alba5}
There exists a map 
$CB(M)\times CB(N)\longrightarrow CB(M\sharp N)$ induced by the 
connected sum of cubulations. 
\end{alba5}
As a consequence $\sharp$ induces also a composition map
$\widetilde{CB}(M)\times \widetilde{CB}(N)\longrightarrow 
\widetilde{CB}(M\sharp N)$.   On the other hand there is a 
natural composition map on the sets of bordisms of immersions, by 
using the connected sum away from the immersions.  We  prove
that the map $I$ is functorial: 
\newtheorem{alba55}[alba1]{Theorem}
\begin{alba55}
We have a commutative diagram 
\[\begin{array}{ccccc} 
\widetilde{CB}(M)      & \times & \widetilde{CB}(N) &\longrightarrow &
\widetilde{CB}(M\sharp N) \\
I\downarrow &  & \downarrow I & & \downarrow I \\
\mathcal I(M) & \times & \mathcal I(N) & \longrightarrow & \mathcal
I(M\sharp N)
\end{array} \]  
\end{alba55}
We believe that the monoid $CB(S^n)$ is actually a group. 
Notice that $\mathcal I(M)$ has a group structure for any $M$ 
(induced from the cohomotopy group structure), but we don't know 
whether this can be lifted to $\widetilde{CB}(M)$.

\section{The proof of theorem 1.4}
Let $K$ be the derivative complex (having the connected
components $K_i$) associated to a cubulation $C$ 
of the manifold $M$ and let  $\varphi:K\longrightarrow M$ be 
 the associated immersion. 
In order to rule out some pathologies we restrict here to the 
{\em combinatorial} cubulations, meaning that
 the star of each vertex (or the link) is a PL ball 
(respectively a PL sphere) of the right dimension.

Let us first show that we have an induced map 
$I:\widetilde{CB}(M)\longrightarrow \mathcal I (M)$. 
Consider the local picture of a bubble move, viewed in the 
boundary of the $(n+1)$-cube. The set of sections  on the boundary 
are intersections of the hyperplane sections of the $(n+1)$-cube with
the faces. Let $B$ and $B'$ be the  $n$-balls  
interchanged by a bubble move. Then the union of 
hyperplane sections in the $(n+1)$-cube yields  
a bordism  between  the  immersions $\varphi_B$ and  $\varphi_{B'}$ 
in the $(n+1)$-ball. Thus immersions associated to equivalent cubulations are 
cobordant.

The immersion $\varphi$ is {\em pseudo-spine} if the
closures of the connected components of the complementary 
$M-Im(\varphi)$ are  balls, and the image of $\varphi$ 
is connected.  The
immersion is called {\em admissible} if each connected
component $L$ of the set  $X_k(M,\varphi)$ of $k$-tuple points 
is a PL ball and  we have 
$cl(L)\cap X_{k+1}(M,\varphi)\neq \emptyset$
(for all  $k\leq n-1$, where $cl$ denotes the closure).  

Observe that any nc-immersion is cobordant to 
an admissible pseudo-spine immersion.  
In fact consider a set of small (bounding) spheres embedded in $M$, 
transverse to 
the immersion $\varphi$. If they are sufficiently small they 
cut the connected components of the complementary $M-Im(\varphi)$ 
into balls. We add sufficiently many so that all 
connected components of the strata $X_k(M,\varphi)$ are divided into
balls by the additional spheres. 
The immersion $\varphi'$ whose image consists
of the union of $Im(\varphi)$ with the spheres is  
cobordant to $\varphi$. In fact let us choose  
some balls in $M\times [0,1]$ bounded by the spheres  
in $M=M\times \{1\}$. The required cobordism is obtained by 
putting the balls in standard position with respect to 
$Im(\varphi)\times [0,1]$. 

Furthermore we want to associate  a cubulation $C$ to the immersion $\varphi$ 
such that $I(C)=\varphi$. 
For an admissible pseudo-spine immersion one takes the corresponding
cell decomposition of the manifold, then the dual decomposition is
just the cubulation we are looking for. The pseudo-spine condition 
implies that the dual decomposition has the structure of a cubical complex
and the admissibility is required for the cubulation comes 
endowed with a natural marking $\mid C\mid \to M$.

\section{Embeddable and mappable cubulations}
\subsection{Mappable cubulations are equivalent}
The proof of  Theorem 1.5 goes as follows: 
we show first that embeddable cubulations are
equivalent and then that a 
mappable cubulation is equivalent to an embeddable one. 

For the first step consider $N$ sufficiently large so that both cubulations $P$
and $Q$ embed into 
the standard cubulation ${\bf R}^N_c$ of the Euclidean space and such 
that there  exists an ambient isotopy carrying $P$ into $Q$.  The
image of $P$ during the isotopy is denoted by $P_t$. 
We define an approximation $P_t^{st}$ of the manifold $P_t$, which is  
a cubulated sub-manifold of ${\bf R}^N_c$, 
it is sufficiently close to $P_t$, and when $t$ varies the 
family $P_t^{st}$ is either locally constant or changes around a 
"critical value" by a bubble move. 
We realize an arbitrarily fine  approximation 
by taking the  cubical structure be  
${\bf R}^N_c[\varepsilon]$, based on cubes whose edges 
are of length $\varepsilon$, for small $\varepsilon$. 
A way to do that is to divide each cube of the 
lattice into $2^N$ equal cubes. 
Then the  initial cubulations $P$ and $Q$ are replaced by 
some iterated doublings, say $2^mP$ and $2^mQ$.   It remains 
to prove  that $2^mP$ is equivalent to $P$. 
\newtheorem{alb}{Proposition}[section]
\begin{alb}\label{gigi}
Let $P$ and $Q$ be two cubulations of a PL manifold which are 
embedded in the standard cubical lattice ${\bf R}^N_c$.
Then there exists $m$  arbitrarily large,  such that 
$2^mP$ and $2^mQ$ are bubble equivalent. 
\end{alb}
\newtheorem{alb2}[alb]{Proposition}
\begin{alb2}
If $P$ is embeddable, then for big enough $m$  the iterated 
doubling $2^mP$ is bubble equivalent to $P$. 
\end{alb2}
\noindent Notice that the analogous statement for np-bubble moves is
false in general. 
\newtheorem{albic}[alb]{Proposition}
\begin{albic}
Let $P$ be a mappable cubulation. Then $P$ is np-bubble equivalent to
an embeddable cubulation.
\end{albic}

\subsection{The proof of Proposition \ref{gigi}}
Consider a cubulation $X\subset {\bf R}^N_c[\varepsilon]$, of
codimension at least 1. Let us denote by $\Lambda[\varepsilon]$
 the union of all hyperplanes defining the cubulation 
${\bf R}^N_c[\varepsilon]$, which  can be 
written as a disjoint union of the  strata $\Lambda[\varepsilon]^{(m)}$
consisting of   all open codimension $m$ cubes.   
Set $C$ for the cube given by the equations 
$\{ |z_j|\leq 1, j=1,N\}$. 
Let $\overline{C}$ be the $(N-1)$-complex obtained from $\partial C$ by adding 
the hyperplanes $\{z_j=0\}$. 
Denote by  $W$ the star of the origin in $\overline{C}$ i.e. the union 
of cells having the form 
$W_{k,\mu}=\{ z_k=0,  \mu_j z_j\geq 0, |z_j|\leq 1, \forall j \}$, where
$\mu_j\in\{-1,1\}, \; \forall j$.
\newtheorem{halbicr}[alb]{Definition}
\begin{halbicr}
The disk $D$ is a standard model in  $C$ if $D$ is properly embedded 
in $C$ (and transverse to $\partial C$), $D$ is contained in $W$ and 
the origin lies in $int(D)$.  
\end{halbicr}   

\newtheorem{halbic}[alb]{Definition}
\begin{halbic}
Let $C$ be a N-cube of the cubulation  ${\bf R}^N_c[\varepsilon]$. We
say that $X$ is standard with respect to C if the following conditions
are fulfilled: 
\begin{enumerate}
\item $X$ is transversal to $\Lambda[\varepsilon]\cap C$. 
\item  There exists an isotopy supported on 
$int(C)\cup (\Lambda[\varepsilon]^{(0)}\cap C)$, if the codimension 
of $X$ is at least 2, and respectively 
$int(C)\cup (\Lambda[\varepsilon]^{(0)}\cup
\Lambda[\varepsilon]^{(1)}\cap C)$,
if the codimension is precisely 1, which transforms $X\cap C$ 
in a standard model. 
\end{enumerate}
Finally $X$ is standard (or in standard position) with respect to 
${\bf R}^N_c[\varepsilon]$ 
(or $\Lambda[\varepsilon]$) if $X$ is standard w.r.t. all cubes. 
\end{halbic}

Observe first that a plane $L$ transverse to the boundary of $C$ is
in standard position w.r.t. $C$.
In fact using a recurrence argument 
the intersection of $L$ with any face $F\subset \partial C$ 
can be put in standard position by means of an isotopy. 
The union of  standard models for $L\cap F$ over the faces is 
 a PL sphere and the cone centered at the origin on it is 
isotopic to $L\cap C$ and hence it is a standard
model. 

Consider  a submanifold  $X$ which is standard w.r.t. the 
lattice $\Lambda$. Then 
there is an isotopy  transforming  $X$ into $X^{st}(\Lambda)$, where 
 $X^{st}(\Lambda)$ intersects each cube along a standard model.
 In fact this can be done in each cube that $X$ intersects  
and the standard models for different cubes have
 disjoint interiors. It suffices to check  the compatibility 
of the boundary gluings: if $X$ cuts two adjacent cubes $C$ and $C'$ 
then the standard models of $X\cap C$ and $X\cap C'$ can be glued
together. First the  neighborhood of the common face 
is determined by the standard model of $X$ inside the
face. Hence we have to prove the uniqueness of 
the  standard model for $X\cap C\cap C'$, which  follows by a 
recurrence argument on the dimension. 
Therefore the cubical complex $X^{st}(\Lambda)$
is uniquely determined by $X$ and $\Lambda$, and it will be called
{\em the standard model} of $X$ w.r.t. $\Lambda$. 

Let  consider the cubulations $P$ and $Q$ of a  
$n$-manifold $M$, embedded in ${\bf R}^N_c$. 
There exists an isotopy with compact support  
$\varphi:{\bf R}^N\times[0,1]\longrightarrow {\bf R}^N$, 
with $\varphi_0$ being the identity and $\varphi_1(P)=Q$. 
For big enough $N$ one can choose the isotopy
$\varphi\mid_{P\times[0,1]}$ 
to be an embedding. 
We assume that the cubulations $P$ and $Q$ are embedded in the
standard cubulation given by the affine lattice $\tilde\Lambda=\Lambda
+(\frac{1}{2},\frac{1}{2},..., \frac{1}{2})$, which has the origin 
translated  into $(\frac{1}{2},\frac{1}{2},..., \frac{1}{2})$. 
Notice that a cubulation $P\subset\tilde\Lambda$ is automatically 
 in standard position 
w.r.t. $\Lambda$. Moreover the intersection of $P$ with each cube 
of $\Lambda$ is a standard model and so $P^{st}(\Lambda)=P$. 
We obtained a submanifold $Y=\varphi(P\times [0,1])$ 
of ${\bf R}^N$, whose boundary $\partial Y$ is in standard 
position w.r.t. the lattice $\Lambda$. 
We claim that for big enough  $m$ there exists  an isotopy 
carrying $Y$ into $X$, such that   $\partial X$ and $\partial Y$  are 
isometric and $X$ is  in standard position w.r.t. the lattice 
$\Lambda[2^{-m}]$. 
There exists a subdivision 
$\Lambda[2^{-m}]$ (for large $m$), such that $Y$ becomes standard w.r.t.
$\Lambda[2^{-m}]$, after a small isotopy which is identity near the 
boundary. 
We translate  $Y$ into the lattice whose origin is at
$(2^{-m-1},2^{-m-1}, ..., 2^{-m-1})$. Then $P$ and $Q$ transform into  
 $2^mP$ and $2^m Q$, and they  are  in standard
position w.r.t. $\Lambda[2^{-m}]$.  The last condition is an open condition, 
so we can keep fixed $Y$ near the boundary during the isotopy. 
Let us denote by $Z$ the tube describing an 
isotopy  between $2^mP$ and $2^mQ$, which is 
in standard position w.r.t. $\Lambda[2^{-m}]$. There exist
topologically trivial tubes $Z_1$ between $P$ and $2^mP$, and respectively 
$Z_2$ between $Q$ and $2^mQ$, which are in a standard positions
w.r.t. $\Lambda[2^{-m}]$. 
Then set  $X=Z_1\cup Z\cup Z_2$. 

Therefore we derived a PL cylinder $X\subset \Lambda[2^{-m}]$ which
interpolates between $P$ and $Q$. 
Notice that the cubical structure of $P$ in $\Lambda[2^{-m}]$ is 
that of $2^mP$ in $\Lambda$. 
In general one cannot  shell the
boundary from $P$ to $Q$. 
The tube $X$ carries a PL foliation by 
submanifolds $P_t =\varphi_t(P)$, where $\varphi$ states for the
isotopy carried by $X$.
The leaf  $P_t$ does not contain {\em flat directions} in the 
cube $U$ if $P_t$ does not contain 
any  segment parallel to some vector in $\partial U$. 
Using a small isotopy one can get rid of flat directions in  
all leaves  $P_t$.
Furthermore there exists some $m$ 
such that either $P_t$ is standard w.r.t. $U\cap
\Lambda[2^{-m}]$, or else $P_t$ contains vertices of the lattice 
$\Lambda[2^{-m}]$. 
The first alternative would  hold 
if we are allowed to move slightly $P_t$, using an arbitrary small 
isotopy (in order to achieve the transversality).
On the other hand the leaf is not transverse iff  
it contains  vertices from $\Lambda[2^{-m}]$, because 
there are not flat directions. 
Further  ``to be in standard position'' is an open condition 
and so one can choose the constant $m$ such that for each $t$ 
either the leaf $P_t$ is in standard position w.r.t. the cube 
$U \cap \Lambda[2^{-m}]$, or else $P_t$ contains vertices from 
$U \cap \Lambda[2^{-m}]$. The set of those exceptional $t$  
for which the second alternative holds is finite because  
each critical leaf contains at
least one vertex from $U \cap \Lambda[2^{-m}]$ and the different leaves
are disjoint. 
Observe  that we can change the isotopy $P_t$ such that no
exceptional leaf contains
more than one vertex from  $\Lambda[2^{-m}]$, while keeping all the
other properties we obtained upon now. Furthermore there is such a $m$
which is convenient for all cubes that $X$ intersects. 
If one replaces a leaf $P_t$ by the standard model
$P^{st}(\Lambda[2^{-m}])$  
then we get a family of 
cubulations embedded in  ${\bf R}^N_c[2^{-m}]$. This family should be locally
constant, until $t$ reaches an exceptional  value $t_0$ (where the standard
model cannot be defined).  
Set $U$ for the cube of size $2^{-m+1}$ centered at  
the (exceptional) vertex. The intermediary set
$P_{[t_0-\varepsilon, t_0+\varepsilon]}\cap U$, for  
$\varepsilon << 2^{-m}$, is 
a trivial cobordism properly embedded in $U$. 
Using a small isotopy one can change $P_{t_0}\cap U$ into a union of planes 
passing through the vertex.  
If  $dir$ is the set of the $2n$ directions of the coordinate  axes around the
vertex lying in $X$ we put  
$dir_{-}= \{ x\in dir; P_{t_0-\varepsilon} \cap \; x \neq
\emptyset\}, \;$,  
$dir_{+}= \{ x\in dir; P_{t_0+\varepsilon} \cap \; x \neq \emptyset\}$.   
Each such direction is dual to a face of a cube $c$ in the dual lattice 
where the standard models $P_{t}^{st}(\Lambda[2^{-m}])$ live. 
We have 
$dir_- \cup dir_+ = dir$, and $dir_-\cap dir_+=\emptyset$ since 
$P_{t_0}$ separates the directions  
which are cut by  $P_{t_0-\varepsilon}$ from the 
directions cut by $P_{t_0+\varepsilon}$. 
Let $f_-$ and $f_+$ be respectively the union of faces of $c$ 
duals to the directions in $dir_-$ and $dir_+$ respectively. 
Then  $f_-$ and $f_+$ are  PL balls because 
 $P_{t_0-\varepsilon}\cap U$ is a ball, as well as 
$P_{t_0+\varepsilon}\cap U$. We need only to see that both are
non-void. If $f_-$ is empty then $P_{t_0-\varepsilon}$ would 
be contained in a half cube $U_0\subset U$ of the lattice 
$\Lambda[2^{-m}]$. Then 
$P_{t_0+\varepsilon}\cap U_0$ will be  a cylinder, and thus it cannot be 
standard, contradicting our hypothesis.
\begin{figure}
\mbox{\hspace{1cm}}\psfig{figure=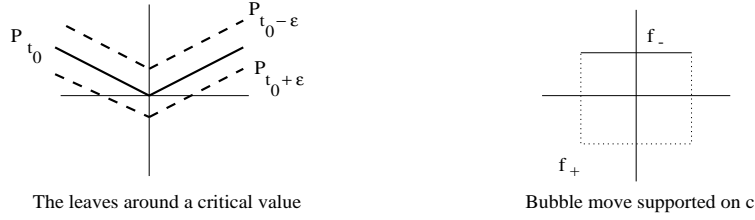,width=10cm}
\caption{The jump of a standard model at a critical value}
\end{figure}
Therefore the standard model $P_{t_0+\varepsilon}^{st}(\Lambda[2^{-m}])$ is obtained
from $P_{t_0-\varepsilon}^{st}(\Lambda[2^{-m}])$ by means of the 
bubble move $f_-\rightarrow f_+$ having the support on $c$. 
This proves Proposition 3.1.

\subsection{The proof of Proposition 3.2}
Actually  a stronger statement concerning 
sub-complexes of the standard lattice ${\bf R}^N_c$ is true. 
We will ask also the bubble moves which pass from one cubulation to the
other to be {\em embedded} in ${\bf R}^N_c$. This means that each 
bubble move which exchanges the balls $B$ and $B'$ has the 
property that the cube bounded by $B\cup B'$ is contained in the 
skeleton of  ${\bf R}^N_c$. 

Two isotopic lattice knots (or graphs) in ${\bf R}^3$ are bubble
equivalent, by means of bubble moves which can be realized on the lattice of 
${\bf R}^3$ and which avoid self-crossings.
Consider  now a lattice $d$-manifold (or complex) $M\subset {\bf R}_c^N$, and a
preferred coordinate axis  defining a height function
$h:{\bf R}^N\longrightarrow {\bf R}$. 
The preimage of the open interval $h^{-1}((n,n+1))$, for $n\in {\bf
  Z}$ is an open PL cylinder because  
an open $n$-cell $e\subset
h^{-1}((n,n+1))$ must be vertical with respect to $h$.
This means that $e= f\times (n,n+1)$, where $f$ is a $(n-1)$-cell whose
projection $h(f)$ is a single point. A horizontal cell is one whose 
image under $h$ is a point. 
Therefore we derive the sub-complexes $A_n^+\subset h^{-1}(n)$, 
$A_n^-\subset h^{-1}(n)$, with $A_n^+\cong A_{n+1}^-$, such that 
$A_n^+\cup A_n^-\cup \{ \mbox{horizontal cells} \} = h^{-1}(n)$ and 
\[cl(h^{-1}((n,n+1))) \cap (h^{-1}(n)\cup h^{-1}(n+1))=
A_n^+\cup A_n^-\cup \{ \mbox{vertical cells} \}.\]
Let $H(n)$ be the union of interiors of the horizontal cells
in $h^{-1}(n)$. Then one can decompose the sets $A_n$ as follows:
$A_n^+=(h^{-1}(n) - H(n)) \cup Z_n$,
$A_n^-=(h^{-1}(n) - H(n)) \cup CZ_n$, 
where $Z_n\cap CZ_n=\emptyset$, and $Z_n\cup CZ_n=\partial (cl(H(n))$.
Using a recurrence argument we  assume that 
$A_n^-$ and $2^kA_n^-$ are bubble equivalent for some 
$k$, by means of the sequence of 
embedded  bubble moves  $X_i$. Then the same
sequences of bubble moves transforms $A_n^+$ into $2^kA_n^+$. 
Since $A_{n+1}^-\cong A_n^+$, there exist some
cone constructions over the  bubble moves $X_i$, which are realized in the 
$(N+1)$-dimensional lattice, and relate  
$A_n^+\times [0,1]$ to  $2^k A_n^+\times [0,1]$, as follows. 
If the bubble move
$X_i$ touches only $A_n^+-Z_n$, then consider the usual cone of $X_i$,
which is also a bubble move in one more dimension. If the bubble move
$X_i$ touches $Z_n$ (or $CZ_n$, on the other side)
then construct an extension with one more dimension for $X_i$, by using the 
horizontal flat. 
The slices  $A_n^+\times [0,1]$, and $h^{-1}((n,n+1))$ can be
glued now back, and we obtain a bubble equivalence between $M$ and
a  $(2^k,2^k,...,2^k,1)$-dilatation, meaning 
that the dilatation acts trivially in the  direction of the chosen axis. 
The same procedure works for any other  coordinate axis. Then the
product of all such dilatations is a homothety of factor 
a  power of $2$, hence the claim follows. 

\subsection{The proof of Proposition 3.3}
A mappable cubulation is not embeddable for two reasons: either it
is non-standard or else the map to the lattice is not
injective. Both accidents can be resolved using np-bubble moves. 
The cubulation  $C$ is $k$-{\it standard} if its $k$-skeleton is
standard, 
i.e. two cubes of dimensions at most $k$ are either disjoint
or else they have exactly one common face.
Assume for the moment that the $n$-dimensional cubulation  
$C$ is $(n-1)$-standard.   
Let $f: C\longrightarrow {\bf R}^N_c$ be a combinatorial map, 
locally an isometry  on each cube. The singularities of the map 
$f$ are therefore either {\em foldings}
or   {\em double points}.  A {\em double point singularity} 
is when  two disjoint cubes  $x$ and $y$ have the same image $f(x)=f(y)$.
The codimension of the cubes is called {\em the defect} of the double points. 
If the defect is positive then $x\subset u \cup v$, $y\subset u'\cup
v'$, where $u$ and $v$ (respectively $u'$ and $v'$) are top
dimensional cubes having a common face. 
A {\em folding} corresponds to a pair of cubes having a
common face and the same image under $f$. 
In order to get rid of double points of positive defect one performs a 
 $b_1$ followed by a $b_2$ move on $u$ and $v$ respectively, in 
an additional dimension. 
If the defect of the double point is zero then we can settle 
using a cubical Whitney trick, as follows. Consider 
$C'$ be obtained by a $b_1$ move on $y$. 
Perform  a 
$b_1$ move on $f(y)$  in  an additional dimension. 
Then   
$D'= b_1(f(C))$ embeds into ${\bf R}^{N+1}_c$, since the  
interiors of the new cubes do not intersect the  cubes 
of $f(C)$. There exists an extension 
$f':C'\longrightarrow D'$ of $f\mid_{C-\{y\}}$ 
sending  $b_1(y)$ onto  $b_1(f(y))$, which  is combinatorial and
has less singularities than $f$. Iterate  this procedure until 
an embeddable  cubulation is obtained.  
In order to solve a folding of two cubes one 
uses a $b_1$ move over one folding cube, and on its image. 
The folding is replaced then  by a double point singularity.

It remains to see how we can use np-bubble moves in order to assume the
$(n-1)$-standardness. Using $(n-1)$-dimensional np-bubble
moves $b_i(n-1)$ one can transform $ske^{n-1}(M)$
into a standard complex. It suffices  to observe that the action of
the $n$-dimensional np-bubble move $b_i(n)$ in one more dimension
agrees with that of $b_i(n-1)$ on the $(n-1)$-skeleton. 
This ends the proof of the Theorem 1.5.    
\newtheorem{albixc}[alb]{Corollary}
\begin{albixc}
The simple cubulations of a manifold $M$ satisfying 
 $H^1(M, {\bf Z}/2{\bf Z})=0$ are bubble equivalent. 
\end{albixc}

\subsection{Np-bubble equivalence and  mappability}
We want to find out whether  a cubulation which is bubble equivalent
to a mappable one is mappable itself. 
In general the answer is negative, and we have to restrict to 
np-bubble moves.  The set of embeddable 
cubulation is not stable under  np-bubbles either. In fact 
an embeddable cubulation may become non-standard,
after performing some $b_k^{-1}$ moves. We will show that this is
the only accident which can happen.
More precisely, we say that $M$ and $N$ are 
{\em standard np-bubble} equivalent if there exists a chain of np-bubble
moves joining $M$ and $N$ among standard cubulations. 
Also the simplicity is preserved by all np-bubble moves but $b_2$, hence
the mappability cannot be np-bubble invariant. 
Two cubulations are {\em simply np-bubble} equivalent if they are
 np-bubble equivalent among simple cubulations. 
Let $M$ be a mappable cubulation and $X\subset M$ be the support of a
$b_2$ move. The  move is 
{\em rigid} with respect to $M$ if there is a combinatorial map 
$f:M\longrightarrow {\bf R}^N_c$ for which $f(X)$ is the union of two 
orthogonal $n$-cubes. Otherwise, either   
$f(X)$ consists of a single cube ($f$ is a {\em folding}) or $f(X)$ is the 
union of two cubes lying in the same $n$-plane ($f$ is {\em flexible}).   
The following result is a coarse converse of the Theorem 1.5:
\newtheorem{albi}[alb]{Proposition}
\begin{albi}\label{albi}
\begin{enumerate}
\item Consider $M$ be mappable and $N$ be a cubulation which is np-bubble 
equivalent to $M$ using only rigid  $b_2$ moves.  
Then $N$ is mappable.
\item If $M$ is embeddable then $b_2(M)$ is mappable.  
\item If $M$ is embeddable and $N$ is standard np-bubble equivalent to $M$ then
$N$ is embeddable.
\item The class of mappable cubulations is closed to simple np-bubble
  equivalences.
\end{enumerate}
\end{albi}
\begin{proof}
Let   $f:M\longrightarrow {\bf R}^N_c$ be a combinatorial  
map into the standard cubulation and  $D_k^+\subset M$ be the support of
a $b_k$ move (i.e. the union of  $k$-cubes).  
Since $f$ is non-degenerate on each cube either 
$f$ is an embedding on $D_k^+$  or else  $f$ is a folding and $k=2$.  
Also  $f$ is always an embedding on the support $D_k^{-}$ of $b_k^{-1}$. 
If the move is rigid (for $k=2$) or $f$ is an embedding on the support  
 then there exists a cube  $C\subset {\bf R}^{N+1}_c$, such that 
$C \cap {\bf R}^N_c = f(D_k)$. 
In fact there is an unique embedding 
of  $D_k^{+}$ ($k\neq 2$) and respectively $D_k^{-}$ into  ${\bf R}^{N+1}_c$, 
up to isometry. 
The map $f$ extends then to a map $\tilde f$ over 
$b_k(M)$ using the cube $C$. 
This ends the proof of the first part of the Proposition. 
\end{proof}

Let us  introduce   some 
notations and definitions from \cite{dss}, for the sake of
completeness. 
Two edges $e$ and $e'$  of a cubic complex $Q$ are said to be 
{\em  equivalent} if there exists a sequence of edges joining them, 
in which any two successive edges 
$e_i, e_{i+1}$ are opposite sides of some square in $Q$.   
An edge equivalence class is
called {\em simple} if all the edges in it belonging to a single cube 
in $Q$ are parallel.  An equivalence
class is called {\em orientable} if all the edges in it can be
oriented such that whenever two equivalent edges are parallel to each
other, their orientations are parallel. 
We consider a  partition $F$  of the edge equivalence
classes  into certain families of classes $F_1,..., F_N$ such that 
each class has a fixed orientation and two 
perpendicular edges are members of equivalence classes that belong to
different families. 
Let $\gamma= (e_1, e_2,..., e_k)$ be an oriented edge path
in the cubulation $Q$. Let $sgn(e_i)$ be +1 if the direction of travel
of $\gamma$ on the edge $e_i$ coincides with the orientation of
$e_i$, and -1 otherwise. Consider then the following 
 formal sum, taking values in the  free ${\bf Z}$-module 
generated by the symbols $F_j$:
\[ D_F(\gamma) = \sum_{i=1}^{k} sgn(e_i) F(e_i)\in {\bf Z}<F_1,F_2,...,F_N> \]
where $F(e)$ is the family to which the edge $e$ belongs.  
We can state now (see \cite{dss}, p. 305-306): 
 \newtheorem{albri}[alb]{Proposition}
\begin{albri}
\begin{enumerate}
\item A simple cubulation  $Q$ maps into the  skeleton of 
${\bf R}^N_c$, where $N$ is the affine dimension of the image, if and only if there
exist orientations of the equivalence classes and a partition $F$ of these
classes into $N$ families such that for  any closed edge path $\gamma$
in $Q$, we have $D_F(\gamma) = 0$. 
\item A simple and  standard cubulation  $Q$ embeds into the  skeleton of 
${\bf R}^N_c$, where $N$ is the affine dimension of the image, if and only if there
exist orientations of the equivalence classes and a partition $F$ of these
classes into $N$ families such that, for an edge path $\gamma$
in $Q$, we have $D_F(\gamma) = 0$ if and only if the path $\gamma$ is 
closed. 
\end{enumerate} 
\end{albri}

We are ready  to prove now the second part:
\begin{proof}
According to  Proposition 3.8 the cubulation $M$ (which is mappable,
hence simple) admits a partition $F$ of the 
equivalence classes of edges such that the development map $D_F$
vanishes on all closed curves. 
The orientations of the edge equivalence classes of $M$ 
naturally induce orientations for $b_2(M)$. 
Further any loop in $b_2(M)$ can be deformed to a loop in $M$ 
hence $b_2(M)$ is  mappable from the  previous
criterion, provided that $b_2(M)$ is simple.
The simplicity is not preserved by $b_2$, in general, but 
 we asked the cubulation  be an embeddable. 
If  $b_2(M)$ won't be simple, then the support of the $b_2$ move should
consist of two {\em twin}
cells. This means that  we have two
cells $e$ and $e'$ which have a common face $f$, and the two  layers
parallel to $f$ containing $e$ and $e'$ coincide. 
But  the layers of an embeddable cubulation  cannot be self-tangent  (see
\cite{dss3}), hence $b_2(M)$ is simple. 
\end{proof}

The third and fourth  statements in Proposition \ref{albi} follow the
same way.

Remark that the collection of non-trivial  homotopy class of immersions 
$\varphi(K_i)$ (of the connected components $K_i$ of the 
derivative complex $K$) is invariant to np-bubble moves.
In particular  $CBB(M)$ is infinite if the manifold $M$ has non-trivial
topology. If we consider two homothetic 
cubulations $X$ and $\lambda X$, then
in general they cannot be np-bubble equivalent, because any non-trivial 
homotopy class appears $\lambda$ times more in the latter.   
Thus the np-bubble 
equivalence may be interesting only for PL-spheres. We believe that 
any two simple cubulations of the sphere are np-bubble equivalent, as
it happens for $S^2$. 

\section{Multiplicative structures}
\subsection{The composition of cubulations}
Assume the manifolds we consider in the sequel  are connected. 
We wish to prove  that the connected sum of cubulations induces a
monoid  structure on  
$CB(S^n)$. 
We believe that $CBB(S^n)$ inherits also a monoid structure. 
Consider the cubulations $C\in C(M)$ and $D\in C(N)$ of the manifolds
$M$ and $N$, and choose 
two cells $e\subset C$, $e'\subset D$. A cell is always a top
dimensional cube in this 
section. Let $t$ be  a length $l\geq 3$ cubical cylinder made 
up from $e\times [0,l]$ by removing the interiors of $e\times\{0\}$ and 
$e\times\{l\}$. 
We define therefore the map ``connected sum of cubulations''
$c_{e,e',t}: C(M)\times C(N)\longrightarrow C(M\sharp N)$, 
by 
\[ c_{e,e',t}(C,D)= C-int(e) _{\partial e \cong \partial e\times \{0\}\subset \partial t}\cup_{\partial e' \cong \partial e'\times \{1\}\subset \partial t}
D-int(e').\]
with the obvious identifications of the boundaries. 

We make some remarks before we proceed to prove the Theorem 1.6. 
We will freely use during the proof the fact that the
tube $t$ can be changed into another tube while it  remains 
mappable. 
The gluing of the tube $t$ requires some 
self-identification of the boundary $\partial t$. We fix an 
arbitrary identification 
$\partial e\times \{0\}\longrightarrow \partial t$. 
Then  the other boundary of the
tube can be glued to $D-int(e')$ in $2^n$ ways corresponding  to the
elements of the symmetry  group $D_n$ of the $n$-cube. 
The relative
twist $tw\in D_n$ measures the difference between two gluings 
on $e'$. Notice that there is no   
there is no canonical choice in gluing $e'$. 
All tubes of length $l\geq 3$ are bubble
equivalent rel boundary so that their length has not to be specified. 
Then the connected sum cubulation depends on the choice of $e$, $e'$
and of the relative  twist.
\begin{proof}
One can assume that $C$ and $D$ are standard cubulations. 
We want  to define a {\em developing map} 
$S:\{$paths in $D\}\longrightarrow \{$cells in $C\}$.
The map $S$ depends on the particular data $e, e', t, tw$ we
chose. 
A {\em path} in $D$ is a sequence of cells starting at $e'$,  
consecutive cells having a common face. It suffices
to define the map $S$ for the trivial path    
then use a recurrence on the length of the path. 
If the path is trivial, consisting of the cell $e'$, we
define $S(e')=e$. Further choose a cell $f'$  having a common
face $u'$ with  $e'$. Let $u$ be the face of $e$ which is 
opposite to  $u'$  using  the tube $t$. Then the  face $f$ neighboring
$e$ and intersecting it along $u$ is by definition $S(\gamma)$, 
where $\gamma$ is the path $(e',f')$.

Let $t'$ be a tube isometric to $t$ which is glued on $C$ and $D$
along $f$ and $f'$ respectively, such that $cl(t\cup t'-t\cap t')$  
is a cylinder on the union of two neighboring cells. Here $t\cap t'=
e\cap e'\times[0,l]$ is the common face of the tubes. 
This condition determines uniquely the gluing twist $tw'$ of $t'$. 
Let $Y$ denotes the cubulation obtained by gluing 
$cl(t\cup t'- t\cap t')$ to $C$ and $D$, along $e\cup e'$ and $f\cup f'$ 
respectively. Then $Y$ is obtained from  $c_{e,e',t}(C,D)$ by a
sequence of $b_2$ moves (along the tube $t'$) and a final $b_3$ move. 
The same procedure transforms $c_{f,f',t'}(C,D)$ into $Y$, hence 
$c_{e,e',t}(C,D)$ and $c_{f,f',t'}(C,D)$ are equivalent. Using a
recurrence one proves that $c_{e,e',t}(C,D)$   and 
$c_{S(\gamma),\gamma(m),t_{\gamma}}(C,D)$ are equivalent for any path
$\gamma$, where $t_{\gamma}$ is a tube isometric to $t$ whose gluing twist
is that induced by $\gamma$. 
\begin{figure}
\mbox{\hspace{1cm}}\psfig{figure=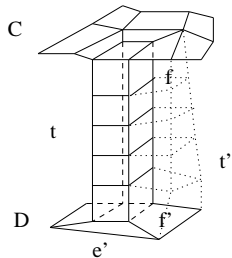,width=3cm}
\caption{The developing map}
\end{figure}

Set 
$ O(e)=\{f' \mbox{ such that } c_{e,e',t}(C,D) \mbox{ is 
  equivalent to } c_{e,f',t'}(C,D) \mbox{ for a suitable twist}\}$. 
We have to prove that,
possibly using bubble moves on the initial cubulations, we have 
$O(e)=C$. 
For any  loop $\gamma$ based at $e$ one knows that $S(\gamma)\in O(e)$. 
Consider a loop  $\gamma$ of length $m$ and assume that the last segment of the curve 
$\gamma\mid_{[m-m_0,m]}$, for $m>m_0\geq 3$ is {\em  straight}.
A curve is {\em straight} if it consists of a sequence of cells $e_i$, 
each cell $e_i$ having a common face $f_i$ with the preceding cell
$e_{i-1}$, so that the faces $f_i$ and $f_{i+1}$ are opposite faces in
$e_i$. 
Then a straight curve $\gamma\mid_{[m-m_0,m]}$ defines a strip
$\Sigma$ consisting of the  cells of the maximal straight extension 
of  $\gamma\mid_{[m-m_0,m]}$. 
Set $f'=S(\gamma)$. The image  of the strip $\Sigma$ under the developing
map is also a strip $\Sigma'\subset C$. 
When a  basepoint cell $e$, and a preferred direction are fixed 
a discrete flow is defined  on the strip $\Sigma$. 
The action of $k\in {\bf Z}$ on the cell $v$ is 
the cell $(k,v)$ which is $k$ steps forward in the given direction,
starting from $v$. Notice that there is also a flow defined  on the
other strip $\Sigma'$.

We claim that  $(2{\bf Z},f')\subset O(e)$ (after using bubble moves).
Let $u=\gamma(r)$, $r\in[m-m_0,m]$, be a cell of $\Sigma\subset D$, 
located between $e$ and $\gamma(m-m_0)$. One performs a $b_1$-move on
$u$ and set   $\gamma'$ for   the natural extension of $\gamma$
to $D'=b_1(D)$.  We can express 
$S(\gamma')$ in terms of the flow on the strip $\Sigma'$ as $(2,
f')$ because  $\gamma'$ is also straight and its
length was increased by  $2$.  Thus $(2{\bf Z}_+,f')\subset O(e)$ holds, 
but  the strip is finite hence  the ${\bf Z}$-action has 
cyclic orbits, implying that $(2{\bf Z},f')=(2{\bf Z}_+,f')\subset
O(e)$. A {\em bicoloring} of $C$  
associates a color $c(e)\in \{0,1\}$ to each cell $e$ such that adjacent cells 
have different colors.  
For any two cells  $f$ and $f'$  there exists a system of  
strips allowing to pass from  $f$ to $f'$ using only the action of
$2{\bf Z}$ on these strips, except for the case when $C$ admits
bicolorings and $c(f)\neq c(f')$. 
However   any  bicolored  cubulation can be transformed into one having no 
bicolorings using  one $b_1$-move. 
This shows that $O(e)=C$. 

Recall that given  a path $\gamma$ between $e$ and $f$ 
one associates  a twist for the corresponding tube over $f$. 
Set $tw_0$ for the twist over $f$ associated to a 
straight path $\gamma=(e,f,g)$. 
Let us perform  a $b_1$ move over $f$. The curve $\gamma$ has several
lifts $\gamma_j$ (not necessarily straight) relating $e$ and $g$
through the cells of $b_1(f)-f$. 
It is simple to check that the set of  twists over $f$ induced by
these paths is  all of $D_n tw_0$. Therefore the connected sum does
not depend on the choice of the twist. This proves the Theorem 1.6.   \end{proof}

\subsection{The compatibility with the bordism composition}
The composition law for immersions is the disjoint union of immersions 
inside the connected sum of manifolds, where the latter is made away from the 
immersions. 
The immersion associated to the connected sum of cubulations 
is obtained from the initial immersions by some surgery which involves
only local data.  
We will prove that this 
surgery can be also realized by a local relative cobordism 
of the associated immersions.  
Consider  the manifolds $M$ and $N$ with their respective 
cubulations $X$ and $Y$.  The connected sum cubulation $X\sharp Y$
is obtained  by removing the cells $e$ from $X$   and $f$ from $Y$. 
One uses $b_1$ moves on  $e$ and $f$ in order  to reduce the
cubulated 
piping tube to a  boundary identification $\partial e=\partial f$. 
The bordism classes of $\varphi_X$ and $\varphi_Y$ are preserved under
these transformations. 
We consider that the connected sum of manifolds is done by means of a 
piping tube which is close to the images 
$Im(\varphi_e)$ and $Im(\varphi_f)$. Both $Im(\varphi_e)$ and
$Im(\varphi_f)$ have as local models the  set of coordinate hyperplanes around 
the origin in ${\bf R}^n$. The surgery which changes 
$I(X\sharp Y)$ into $I(X)\sharp I(Y)$  
excises $Im(\varphi_e)$ and
$Im(\varphi_f)$ and  replaces them by a cylinder with the same
boundary. One can realize the surgery on the piping tube, after pushing
the local models thru it.  A small  isotopy 
moves them outside of a longitude and then the configurations 
embed in the ball obtained from the tube by cutting it along the longitude. 
It suffices  therefore to show that the surgery can be realized by a
cobordism  in the ball which is a product on the boundary. 
For $n=2$ the different slices of the cobordism are given in the 
picture 4. 
\begin{figure}
\mbox{\hspace{2cm}}\psfig{figure=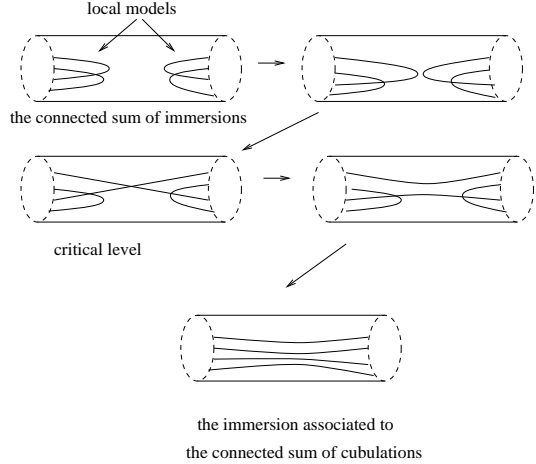,width=7cm}
\caption{A relative cobordism for $n=2$}
\end{figure}
Assume that the local models living in a ball have the corresponding 
hyperplane sections parallel to each other.  
Let us consider two parallel hyperplane  sections $u$ and $u'$ 
of the respective 
local models $Im(\varphi_e)$ and $Im(\varphi_f)$.  
Then one constructs the neighborhood of  $\varphi_{X\sharp Y}$ 
around the piping tube by 
removing the interiors of $u$ and $u'$ and gluing back  the 
cylinder of   boundary $\partial u \cup \partial u'$.
This transformation can be realized also  by a cobordism of immersions
since it is represented by 
the local picture around a critical point of index 1, where the 
images  of the immersions are the non-critical levels before and 
respectively  after 
passing through the critical value.  The local picture  
can be made transversal with respect to the 
other coordinate hyperplanes of the immersions which were left untouched.  
One composes the cobordisms associated to  
all $n$ pairs of parallel hyperplane sections. 
The restrictions on the boundary of these cobordism are products 
and so we can glue the composition with the  
trivial cobordism outside the local pictures. Thus we derived a cobordism between 
$\varphi_{X\sharp Y}$ and $\varphi_{X}\sharp \varphi_Y$.

Notice that, in general, the map
$CBB(M)\times CBB(N)\longrightarrow CBB(M\sharp N)$ should depend on the
length of the tube $t$. The case when $M$ and $N$ are spheres 
could be different however. 

\section{Cubulations of surfaces}
\subsection{The simple cubulations of $S^2$}
Observe first  that the set $CBB(S^2)$ is infinite.
In fact the  number of components with odd
self-intersections of the associated immersion is an invariant, 
because the only move creating new components is $b_1$, and each new 
created component is  an embedded
circle. Set $ns(C)$ for the sum over the various connected components
$K_i$ of the number of self-intersections.   

The similarity with the Reidemeister moves in the plane  suggests that
$CBB(S^2)$ is the set of  framed circles in the plane. 
Unlike the case of Reidemeister moves, the image of the immersion
remains connected, so its components
cannot be  separated using bubble moves. 
On the other hand the move $b_2$ can create/annihilate a pair of
self-intersections. It is not clear that the 
singularities can always be paired such that suitable np-bubble moves 
destroy all pairs of singularities and so 
each transformed circle has  
$ns(K_i)(mod \; 2)\in \{0, 1\}$ singularities. If
this is true it will remain to  prove that all configurations of
circles among which there are exactly  $m$ singular  circles are 
np-bubble equivalent.
This would establish an isomorphism between $CBB(S^2)$ and 
${\bf Z}_+$. 
We are able to prove this statement for the 
case $m=0$ (and $m=1$):
\newtheorem{lolo}{Proposition}[section]
\begin{lolo}
The simple cubulations of $S^2$ are np-bubble equivalent.
\end{lolo}
\begin{proof} 
Call a disk bounded by two disjoint arcs a {\em biangle}. A biangle is
{\em tight} if no other arc intersects its interior. 
Observe  that using np-bubble moves one can transform all minimal 
biangles into tight biangles. We use  isotopies  of the boundary of
the biangle ($b_2$-moves) reducing the number of squares contained in 
the disk  and $b_3$ moves. 

The possible configurations of
tight biangles can be rather complicated, even if the cubulation is 
simple.  Let us show that the tight triangles 
can travel along the edges (sliding). This means that a regular
homotopy moving a  tight biangle along the arcs can be realized  by 
 np-bubble moves. The proof of this claim
is contained in the picture below. We marked in a little rectangle the area on which the bubble move acts:

\vspace{0.5cm}
\mbox{\hspace{4cm}}\psfig{figure=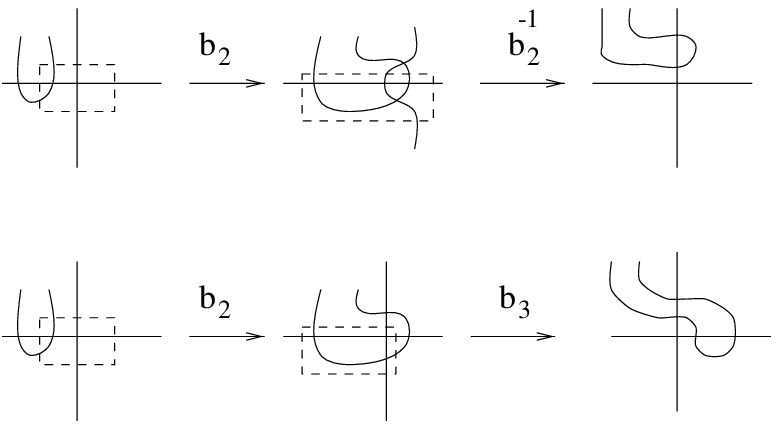,width=7cm}
\vspace{0.5cm}

Set $X$ for the union of two transversal arcs. An arc which intersects
three times one branch of $X$ and then the other arc can be simplified
by a composition of bubble moves which we call a S move:

\vspace{0.5cm}
\mbox{\hspace{4cm}}\psfig{figure=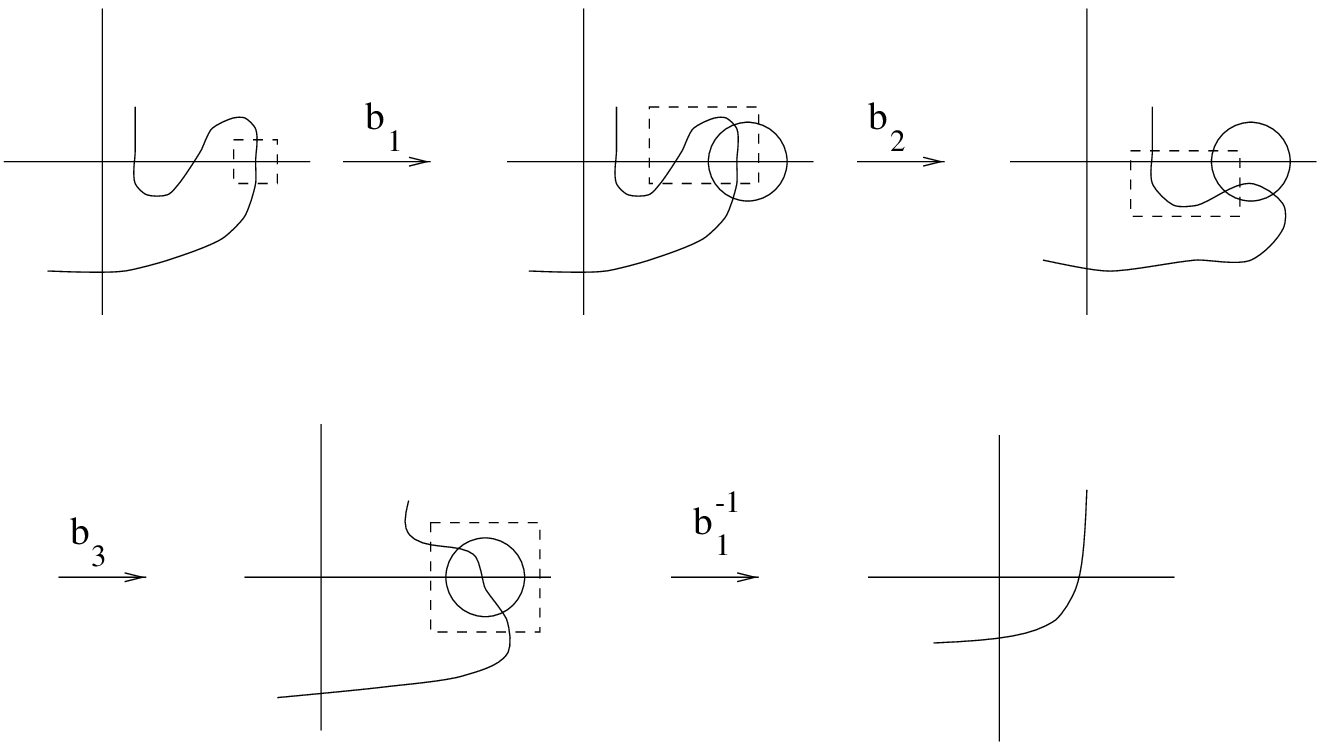,width=8cm}
\vspace{0.5cm}

Now it is easy to obtain that the set of semi-simple cubulations mod np-bubble moves is equal 
to the set of simple cubulations mod np-bubble moves. 
In fact any tight biangle coming from the same component can be transformed
into tight biangles on different components:  

\vspace{0.5cm}
\mbox{\hspace{4cm}}\psfig{figure=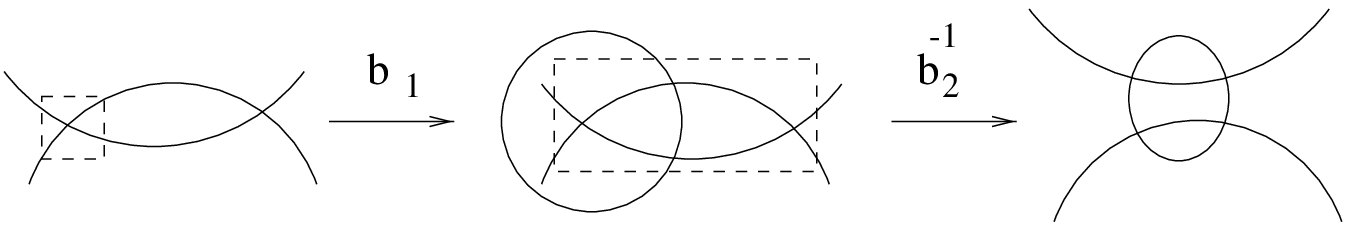,width=8cm}
\vspace{0.5cm}

Assume  that the cubulation is simple. 
The union of the circles $K_i$ for $i\geq 2$,  divides the  sphere
into polygonal faces, each of them  having at least two  vertices. 
We claim that, either there are no biangles involving an arc from
$K_1$ or else $K_1$ is  a small circle contained in the union of
two faces which intersects minimally (i.e. twice) the common edge. 
Indeed, suppose that there exists a minimal tight biangle. 
Consider  the face which is adjacent to the biangle 
and not containing it. If $K_1$ does not satisfy the claim then 
we can simplify the biangle using slidings and S-moves,  
as it is shown below:

\vspace{0.5cm}
\mbox{\hspace{4cm}}\psfig{figure=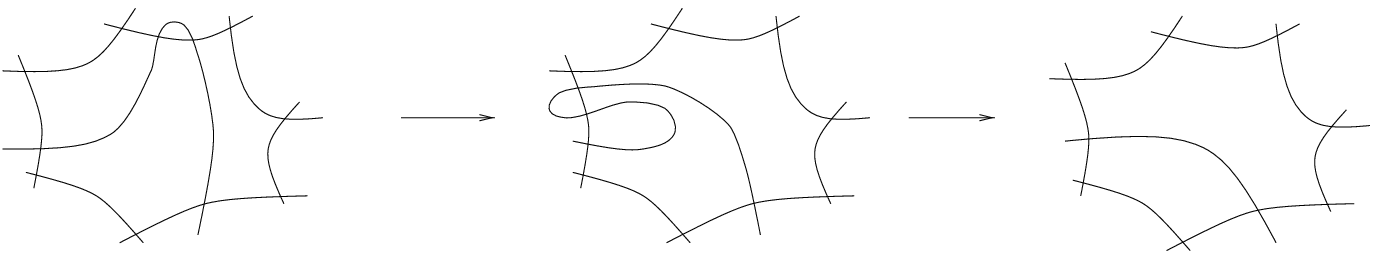,width=8cm}
\vspace{0.5cm} 

Thus $K_1$ can be transformed into a circle which does not 
support any biangle. 
Further use of the moves $b_2$ and $b_3$ allows us to isotope  
$K_1$ into the union of two faces as claimed.  
One continues the simplification  procedure with the other components
$K_i$ and the Proposition 5.1 follows.   
\end{proof}
A similar proof works when  the circle $K_1$ is allowed to have one 
self-crossing. Now $K_1$ is transformed into a 
a figure eight contained in the union of two faces which  
intersects minimally their common edge and forms one biangle. 

\subsection{The cubulation group $CB(S^2)$}
In order  to prove that $CB(S^2)={\bf Z}/2{\bf Z}$
one has  to get rid of  those self-intersections
which are not cancelling pairs. The following figures
describe the simplifications obtained with the additional move
$b_{3,1}$ for two adjacent self-intersections:

\vspace{0.5cm}
\mbox{\hspace{1cm}}\psfig{figure=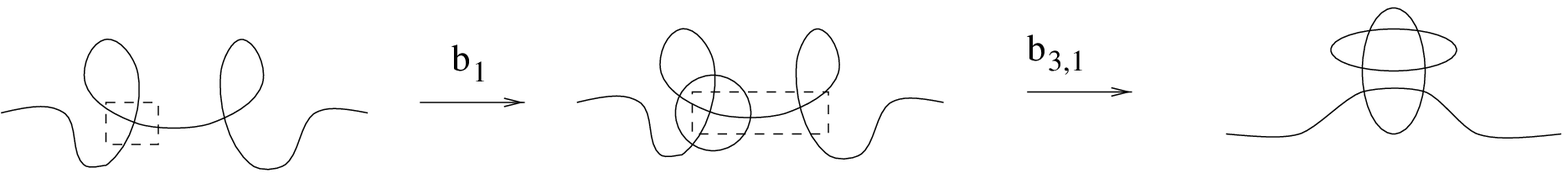,width=12cm}
\vspace{0.5cm}

Further, if the two self-intersections are separated by additional  arcs
then use $b_2$-moves  and slide across these arcs. 
It follows that any dual graph can be transformed using bubble moves 
into one satisfying $ns(K_i)\in \{0,1\}$, for each component $K_i$.

The next step is to show that two components ``can be added''. Let
$K_1$ and $K_2$ be two components having non-void intersection. 
There exists an equivalent configuration  in which $K_1$ and $K_2$ 
are replaced by the circle $K_1+K_2$ verifying 
$ns(K_1+K_2)=ns(K_1)+ns(K_2)$, the 
other $K_i$'s ($i\geq 3$) are left unchanged and the additional 
components have $ns=0$. It suffices to do that in the 
case when  $ns(K_1)=ns(K_2)=1$. Then the kinks can be added and 
transformed in consecutive self-intersections which have already been solved. 

\vspace{0.5cm}
\mbox{\hspace{2cm}}\psfig{figure=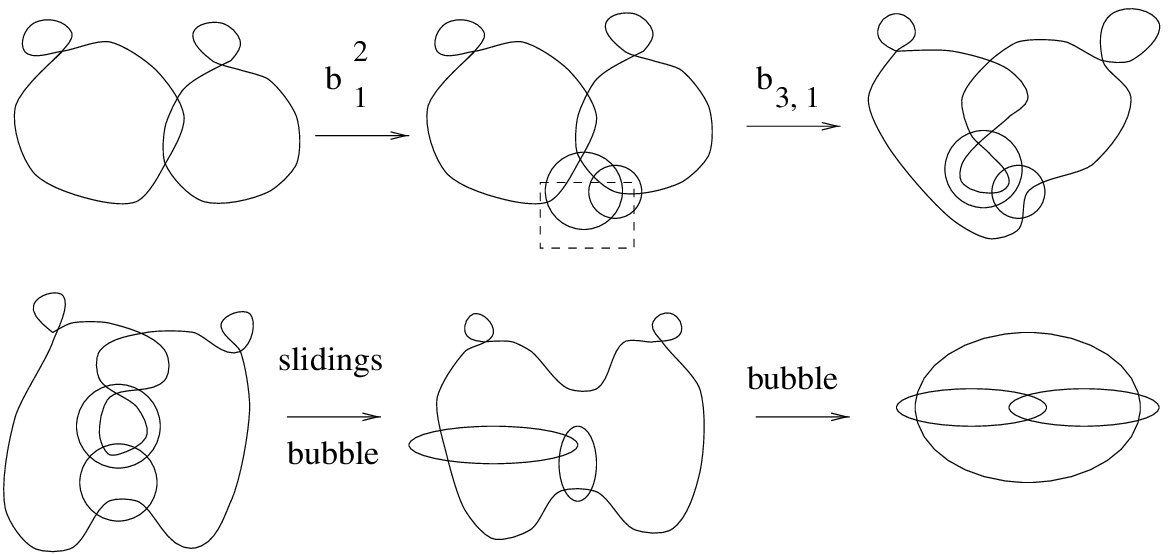,width=10cm}
\vspace{0.5cm}

Thus, if one adds those components which are not embedded, we obtain 
a planar graph whose components are embedded, except possibly for 
one which has $ns(C)\in\{0,1\}$. We have therefore  
$ns(C)=f_0(mod \, 2)$. 
If $f_0=0(mod \, 2)$ then the result of the previous section shows
that all these configurations are np-bubble equivalent. 
If $f_0=1(mod \, 2)$ it means that $C$ is a figure eight in the plane.
The remark after the proof of Proposition 5.1  ends the proof of the 
Theorem 1.3.

\bibliographystyle{plain}

\end{document}